\documentclass[a4paper, pagesize, DIV12
	       ]{scrartcl} 
\usepackage{typearea}

\usepackage[latin1]{inputenc}
\usepackage[ngerman, english]{babel}

\usepackage{amssymb}
\usepackage{amsmath}
\usepackage{amsthm}
\usepackage{mathtools}

\usepackage{graphicx} 

\title{Hoppe trees, random recursive sets and their barycentre}
\author{Mathias Rafler\thanks{rafler@ma.tum.de}\\
}

\renewcommand{\l}{\ell}

\renewcommand{\Pr}{\mathbb P} 
\newcommand{\Ex}{\mathbf E} 
\newcommand{\Var}{\mathbf{Var}} 
\newcommand{\Cov}{\mathbf{Cov}} 

\newcommand{\hoppe}{\mathcal{H}}

\newcommand{\Poi}{\mathbf P\!\!} 
\newcommand{\Nd}{\mathcal N} 

\newcommand{\Esig}{\mathcal{E}}

\newcommand{\M}{\mathcal M}

\newcommand{\N}{\mathbb N}
\newcommand{\Z}{\mathbb Z}
\newcommand{\R}{\mathbb R}

\renewcommand{\d}{\mathrm{d}}

\newcommand{\eqa}[1]{\begin{align*}#1\end{align*}}
\newcommand{\equ}[1]{\begin{equation*}#1\end{equation*}}
\newcommand{\eqan}[1]{\begin{align}#1\end{align}}
\newcommand{\equn}[1]{\begin{equation}#1\end{equation}}

\newcommand{\distreqq}{\stackrel{\mathrm{d}}=}

\theoremstyle{definition}
\newtheorem{defdefinition}{Definition}[section]

\theoremstyle{plain}
\newtheorem{defsatz}[defdefinition]{Theorem}
\newtheorem{defsatzdef}[defdefinition]{Theorem and Definition}
\newtheorem{defprop}[defdefinition]{Proposition}
\newtheorem{deflemma}[defdefinition]{Lemma}
\newtheorem{deffolgerung}[defdefinition]{Corollary}

\theoremstyle{remark}
\newtheorem{defbemerkung}[defdefinition]{Remark}

\newcommand{\satz}[1]{\begin{defsatz}#1\end{defsatz}}
\newcommand{\satzn}[2]{\begin{defsatz}[#1]#2\end{defsatz}}

\newcommand{\prop}[1]{\begin{defprop}#1\end{defprop}}

\newcommand{\bem}[1]{\begin{defbemerkung}#1\end{defbemerkung}}
\newcommand{\lemma}[1]{\begin{deflemma}#1\end{deflemma}}
\newcommand{\lemman}[2]{\begin{deflemma}[#1]#2\end{deflemma}}

\newcommand{\korollarn}[2]{\begin{deffolgerung}[#1]#2\end{deffolgerung}}

\numberwithin{equation}{section}

\begin{document}

\maketitle

%

\begin{abstract}
We consider a recursively defined random set of points and its barycenter, where the random set is constructed by the following inductive rule: Given a realization of $n-1$ points, one of them is picked at random and serves as a source the $n$-th point.  We discuss the asymptotic behaviour of the barycentre of this random set. The main analysis relies on the analsis of Hoppe trees, for which we derive a limit theorem for the joint distribution of total length and Wiener index.
\end{abstract}

\section{The Hoppe construction of a random recursive set}

In~\cite{fH87}, Hoppe constructed a P\'olya-like urn model by introducing a kind of magic ball, which in case of being drawn, forces to add a ball of a new colour to the urn. If a ball of any other colour is drawn, one follows the mechanism of P\'olya's urn. One parameter in this model is the weight $\theta>0$ of the magic ball, which is allowed to differ from the weight of the other balls. In the beginning, a large $\theta$ favours the introduction of new colours, but on the long run this influence looses importance.

If each ball get an unique identifier and one memorizes the parent of each ball, one constructs some kind of ancestral tree, which is called Hoppe tree in~\cite{LN12}. The magic ball becomes the root and every branch of the tree is what in Hoppe's urn model is the colour. For $\theta=1$ this tree model is the random recursive tree.

A similiar construction appears in the context of point processes, a particular example of Papangelou processes: One aims at constructing a point process by specifying simply a conditional intensity which is the intensity of adding further points given a realized point configuration. The basic example is the Poisson process $\Poi_\rho$ with intensity measure $\rho$, for which the conditional intensity is exactly $\rho$, independent of the given point configuration. Moreover, the Poisson process is uniquely specified by this conditional intensity. Considering the point configurations as point measures, Zessin constructs in~\cite{hZ09} a point process with conditional intensity $z(\rho+\mu)$ ($\mu$ is the point measure given by the point configuration including multiplicities), $z\in(0,1)$ a real number, i.e. presence of points yields rewards for the intensity. This point process is unique and called P\'olya sum process.

As Leckey and Neininger~\cite{LN12} added more information to the Hoppe urn, here we add some more information for a spatial component like Zessin did. Instead of just rewarding the points in space where some other points already are, we give additional weight to set close to that point.

Construct a sequence of random variables $(X_n)_{n\geq 0}$ in the Euclidean space $E$ in the following way: Fix some $\theta>0$ and a stochastic kernel $\kappa$ from $E$ to $E$. For $n\in\N$ let $(J_n)_{n\in\N}$ be a sequence of independent random variables with distribution
\equ{
  \Pr(J_n=0)=\frac{\theta}{\theta+n-1},\qquad \Pr(J_n=k)=\frac{1}{\theta+n-1},\qquad k\in\{1,\ldots,n-1\}.
}
We write $J_n\sim \hoppe_n(\theta)$. Note that for $\theta=1$, $\hoppe_n(\theta)$ is the uniform distribution on $\{0,\ldots,n-1\}$. Finally let $X_0=0$ and define recursively
\equ{
  \Pr(X_n\in B|X_0,\ldots X_{n-1},J_1,\ldots, J_n)=\kappa(X_{J_n},B). 
}
In other words: $X_0$ acts as a source or root at $0$, and given $X_0,\ldots, X_{n-1}$, $X_n$ is given by a jump from a randomly picked point among the given ones according to $\kappa$.


There are several immediate examples of jump kernels $\kappa$ one might think of:
\begin{enumerate}
  \item the normal kernel $\kappa(x,\,\cdot\,)=\Nd(x,\sigma^2)$ for some fixed $\sigma^2>0$, $x\in\R$,
  \item a Poisson kernel $\kappa(x,\,\cdot\,)=x+\Poi_\lambda$ for some fixed $\lambda>0$, $x\in\N$,
  \item the shift kernel $\kappa(x,\,\cdot\,)=\delta_{x+1}$, $x\in\N$,
  \item simple random walk kernel $\kappa(x,\,\cdot\,)=\tfrac{1}{2}\delta_{x-1}+\tfrac{1}{2}\delta_{x+1}$, $x\in\Z$.
\end{enumerate}
Particular feature of the first two kernels is that their iterations are of the same structure and the one-dimensional distributions will be mixed normal and mixed Poisson distributed. In the third example a chosen point produces another one directly at the next location on the right. First and last example are centered ones.

We call $\kappa$ \emph{covariant} if
\equ{
  \int f(y)\kappa(x,\d y)= \int f(y+x)\,\kappa(0,\d y),
}
If we further define
\eqa{
  m(x)&=\int y\,\kappa(x,\d y),\\
  v(x)&=\int\bigl(y-m(x)\bigr)^2\,\kappa(x,\d y),
}
then for covariant $\kappa$, modulo existence, by
\equ{
  m(x) = \int y+x\,\kappa(0,\d y) = m(0)+x
}
and
\equ{
  v(x) = \int \bigl(y-m(0)\bigr)^2\,\kappa(0,\d y) = v(0),
}
the expected offspring location $m(x)$ of a point at $x$ is fixed relative to $x$ as well as its variance $v(x)$ is constant in space. In any of the given examples, $\kappa$ is covariant, and moreover in the first and last example $m$ is the identity. If $m(x)=x$, we say that $\kappa$ is \emph{centered}.

Our aim is to discuss the barycentre $S_n$ of the random recursive set $\{X_0,\ldots,X_{n-1}\}$ and its asymptotic behaviour as $n\to\infty$. In particular we show that for the normal kernel as $n\to\infty$, the barycentre is mixed normally distributed with the mixing measure given implicitly and only depending on the tree constructed from $J_1,J_2,\ldots$. The tool is the limit theorem for the joint distribution of the length and the Wiener index of a random recursive tree in~\cite{rN02}, an application of the multivariate contraction method, which yields an implicit representation of the limiting distribution. This also yields a limit theorem for the joint distribution of length and the Wiener index of the Hoppe tree.

The note ist structered as follows: We give a random-walk like representation of the given recursive set in section~\ref{sect:rwrep}, which yields that the joint distribution of the sequence $(X_n)_{n\geq 0}$ is a mixture of distributions, and the mixing measure is given by a Hoppe tree. In particular for the normal kernel $(X_n)_{n\geq 0}$ is a mixed Gaussian process, and therefore its barycentre should, if well-defined, be mixed normal with a random variance or even normal. It turns out that the former case holds, although the variance of the variance is rather small and vanishes as $\theta\to\infty$. In this case formally only $X_0=0$ produces descendents, and therefore for the barycentre the law of large numbers holds. These results are given in~\ref{sect:bary}; since the proofs substantially rely on the analysis of the Hoppe tree, they are postponed to the Hoppe tree section~\ref{sect:limits}.

\section{A random walk-like representation\label{sect:rwrep}}

Firstly we derive a representation of the distribution of any finite sequence $(X_1,\ldots,X_n)$ in terms of independent random variables. Particularly, the distribution of this vector is a random linear transformation of a product measure provided $\kappa$ is covariant. 

The Hoppe tree was introduced in~\cite{LN12} as the tree constructed in the following way: Start with a root which gets a weight $\theta>0$ and add subsequently new vertices, which get weight 1, by choosing a present vertex proprtional to its weight and attach a new vertex to the chosen one. Note that $J_1,\ldots,J_{n-1}$ generates a Hoppe tree of $n$ vertices, which is rooted at 0 and for which $J_k$ is the parent of $k$.

Conditioned on $\mathcal{J}=\sigma\bigl((J_k)_{k\geq 1}\bigr)$, the sequence $(X_n)_{n\in\N}$ can be described very detailed, and many properties just go back to the Hoppe tree.

We need to fix some notation and for that we only consider the sequence $(J_k)_{k\geq 1}$ for the moment. For any $k\in\N$ there exists a unique path $p_k$ from $0$ to $k$, a strictly increasing sequence with its first entry being 0 and its last one being $k$, which one can construct backwards by going from $k$ to its parent $J_k$ and repetition with $k$ replaced by $J_k$. In reverse order this is the path $p_k$. Denote by $D_k=|p_k|$ the length of this path, i.e. the number of steps needed to reach $k$ from the root. For distinct vertices $j$ and $k$ we have $p_j(0)=p_k(0)$ and denote by $\l_{jk}$ the last point of unity of $p_j$ and $p_k$,
\equ{
  \l_{jk}=\sup\{m\in\N:p_j(m)=p_k(m)\},
}
as well as by $D_{jk}$ its distance to the root.

For a vector $(X_1,\ldots,X_n)^t$ write $X^n$.
\satz{ \label{thm:rwrep:fididi} 
Assume that $\kappa$ is covariant, and let $Y_1,\ldots,Y_n$ be iid with distribution $\kappa(0,\,\cdot\,)$. Then 
\equ{
  X^n\distreqq A Y^n,
} 
where the random matrix $A$ has entries $A_{ij}=1_{\{j\in p_i\}}$. Moreover, letting $\mathcal{J}=\sigma(J_k: k\geq 1)$ for $j,k\leq n$,
\equ{
  \Cov(X_j,X_k|\mathcal{J})=D_{jk}\sigma^2,
}
and without condition
\equ{
  \Cov(X_j,X_k)=\Ex[D_{jk}]\sigma^2.
}
}
Thus the random vector $X^n$ can be constructed in two steps: Firstly realize a Hoppe tree of $n+1$ points, secondly attach to each edge a random variable distributed according to $\kappa(0,\,\cdot\,)$ realized independently of each other and independent of the tree, and determine $X_j$ for $j\in\{1,\ldots,\}$ by summing up all of the random variables along the path $p_j$.
\begin{proof}
Conditioned on $\mathcal{J}$, let $k$ be the length of the path $p_j$ from $0$ to $j$, then 
\equ{
  X_j = (X_j-X_{p_j(k-1)})+\ldots+(X_{p_j(1)}-X_0).
}
By construction, the summands are independent and since $\kappa$ is covariant, they are identically distributed. If $(Y_n)_{n\geq 1}$ is a sequence of independent, $\kappa(0,\,\cdot\,)$ distributed random variables, then we can write
\equ{
  X_j \distreqq Y_j+ Y_{p_j(k-1)}+\ldots + Y_{p_j(1)}.
}
Hence we have $X^n= A Y^n$ conditional $\mathcal{J}_n=\sigma\bigl((J_k)_{1\leq k\leq n}\bigr)$ for any $n$ but then immediatly also without the condition.

For the covariance observe that for $j,k\in\N$, the paths $p_j$ and $p_k$ from 0 to $j$ and $k$, respectively, agree up to $\l_{jk}$, i.e. $X_j$ and $X_k$ share the same summands up to $\l_{jk}$, but after that point they are independent. Therefore
\eqa{
  \Cov(X_j,X_k|\mathcal{J})=\Var(X_{\l_{jk}}|\mathcal{J})=|p_{\l_{jk}}|\sigma^2.\qedhere
}
\end{proof}

If we consider the case of a normal jump kernel, then this result implies that the process $(X_n)_{n\in\N}$ is a Gaussian process conditioned on $\mathcal{J}$.

\korollarn{Normal jumps}{ \label{thm:rwrep:gauss} 
Let $\kappa(x,\,\cdot\,)=\Nd(x,\sigma^2)$. Then given $\mathcal{J}$, $(X_n)_{n\in\N}$ is a centered Gaussian process with covariance function $\Sigma_{j,k}=\Cov(X_j,X_k|\mathcal{J})=D_{jk}\sigma^2$.
}

\section{The barycentre of the random recursive set\label{sect:bary}}

The barycentre $S_n$ of $\{X_0,\ldots,X_{n-1}\}$ is the expectation of the empirical distribution of the given set
\equ{
  S_n=\frac{1}{n}\sum_{k=1}^n X_k.
}
Theorem~\ref{thm:rwrep:fididi} implies that $S_n$ has expectation 0 if and only if $\kappa$ is centered and moreover finite variance if and only if $v$ is finite. If moreover $\kappa$ is the normal kernel, then immediatly $S_n$ conditioned on $\mathcal{J}$ is normally distributed by Corollary~\ref{thm:rwrep:gauss}.

For each $n\in\N$ define
\equ{
  T_n = \sum_{k=0}^{n-1} D_k,\qquad R_n=\sum_{0\leq i<j\leq n-1} D_{ij}.
}
$T_n$ is the total length of the tree. By making use of these definitions,
\eqa{
  n^2\Var(S_n|\mathcal{J})&=\Var(X_0+\ldots+X_{n-1}|\mathcal{J})\\ 
   &=\sum_{j=1}^{n-1}\Var(X_j|\mathcal{J})+2\sum_{0\leq i<j\leq n-1} \Cov(X_i,X_j|\mathcal{J})
	=\sigma^2T_n+2\sigma^2R_n.
}
Thus,

\lemman{Conditional distribution of the barycentre for normal jump kernels}{
Let $\kappa(x,\,\cdot\,)=\Nd(x,\sigma^2)$, then the barycenter $S_n$ of $X_0,\ldots,X_{n-1}$ given $\mathcal{J}$ is centered normally distributed for all $n$ with variance
\equ{
  \sigma^2\frac{U_n}{n^2}:=\sigma^2\frac{T_n+2R_n}{n^2}.
}
}

In particular the only influence of the jumps on the variance of $S_n$ is via the variance of the jumps, the remaining part is depends only on the Hoppe tree and is independent of the single jumps. This is not a special feature of normal jump kernel. Denote by $W_n$ the \emph{Wiener index} of the Hoppe tree of size $n$, that is the sum over all distances between all pairs of distinct vertices. Then immediatly follows that $W_n$ is connected to $U_n$ and $R_n$ in the following way
\lemma{ \label{thm:rwrep:rec:Un}
For $n\geq 1$,
\equ{
  U_n = nT_n-W_n,\qquad R_n=\frac{n-1}{2}T_n-\frac{1}{2}W_n
}
where $W_n$ is the Wiener index of the Hoppe tree with $n$ vertices.
}

The recursive structure allows to obtain recursions for the expectations of these random variables. \eqref{en:rek:Tn}~and~\eqref{en:rek:ETn} are already contained in~\cite{LN12}. By $h_n^\theta$ we denote the sum
\equ{
  h_n^\theta = \sum_{j=1}^{n-1}\frac{1}{\theta+j},
}
$\Psi$ is the digamma function. The proof is contained in the next section.
\prop{ \label{thm:barycntr:expect}
Let $\theta>0$.
\begin{enumerate}
  \item The conditional expectations satisfy the recursions, $n\geq 1$,
  	\eqan{
  	  \Ex[T_{n}|\mathcal{J}_{n-1}] &= \tfrac{\theta+n-1}{\theta+n-2}T_{n-1}+1
  	     \label{en:rek:Tn}\\
   	  \Ex[R_{n}|\mathcal{J}_{n-1}] &= \tfrac{\theta+n}{\theta+n-2}R_{n-1} 
   	  		+\tfrac{1}{\theta+n-2}T_{n-1} \label{en:rek:Rn}\\
			\Ex[U_{n}|\mathcal{J}_{n-1}] &= \tfrac{\theta+n}{\theta+n-2}U_{n-1} 
					+\tfrac{1}{\theta+n-2}T_{n-1}+1\\ 
			\Ex[W_{n}|\mathcal{J}_{n-1}] &= \tfrac{\theta+n}{\theta+n-2}W_{n-1}
					+\tfrac{\theta-1}{\theta+n-2}T_{n-1}+n-1 \label{en:rek:Wn}
		}
	\item For $n\geq 1$ the expectaions are
	  \eqan{
  		\Ex T_n &= (\theta+n-1)h_n^\theta=n\log n-\Psi(\theta+1)n+(\theta-1)\log n+O(1)
  			\label{en:rek:ETn}\\
  		\Ex U_n &= (\theta+n)(\theta+n-1)\bigl[\tfrac{2}{1+\theta}-\tfrac{1}{\theta+n-1}
  				-\tfrac{1}{\theta+n}(1+h_{n-1}^\theta)\bigr]\notag\\
  					&=\tfrac{2}{1+\theta}n^2-n\log n+O(n)\\
 			\Ex W_n &= (\theta+n)(\theta+n-1)\bigl[(\theta-1)\bigl(\tfrac{1}{\theta+1}
  				-\tfrac{1}{\theta+n-1}-\tfrac{1}{\theta+n}h_{n-1}^\theta\bigr)+h_n^\theta-1
  				+\tfrac{\theta+1}{\theta+n}\bigr] \notag\\
  				&=n^2\log n
  				-\Bigl[\Psi(\theta+1)+\tfrac{2}{\theta+1}\Bigr]n^2+\theta n\log n+O(n).
  	}
\end{enumerate}
}

The expectations of $nT_n$ and $W_n$ grow both like $n^2\log n$, and even cancel each other such that $\tfrac{1}{n^2}\Ex U_n$ converges to some finite, non-zero limit. It is the consequence of a limit theorem for the joint distribution of $(T_n,W_n)$ that $\tfrac{1}{n^2}U_n$ converges in distribution together with its second moments -- the limit is given implicitely. These limit theorems together with the proofs of the following corollaries are given in the next section. For the random recursive tree we have

\korollarn{Asymptotic behaviour in recursive tree}{ \label{thm:barycntr:Urek}
Let $\theta=1$. Then the sequence $\bigl(\tfrac{U_n}{n^2}\bigr)_{n\geq 1}$ converges in distribution to a random variable $U$, which is a solution of the fixpoint equation
\equn{ \label{eq:barycntr:Urek}
  U\distreqq V^2 U^\ast+(1-V)^2U+V^2,
}
where $U$, $U^\ast$, $V$ are independent, $V\sim \mathcal{U}[0,1]$ and $U\distreqq U^\ast$. The first two moments are
\equ{
  \Ex U=1,\qquad \Ex U^2=\frac{11}{9}.
}
In particular, if $\kappa$ is the normal kernel for $\sigma^2>0$, the barycentre has asymptotically a mixed centered normal distribution with random variance, whose expectation is 1 and variance is $\tfrac{2}{9}\sigma^2$.
}

It is the nice property that if the first branch of a random recursive tree is cut of, then the two subtrees are independent random recursive trees and given their total number of vertices, the number of vertices of the subtrees is uniformly distributed. The analogue result is holds for the case $\theta\neq 1$ with the difference that one of the trees is a random recursive tree und that the size of one subtree given the total size is no longer a uniform distribution. To distinguish between the two trees, we write $U_n'$, $T_n'$, etc. whenever the underlying tree is the Hoppe tree.

\korollarn{Asymptotic behaviour in Hoppe tree}{ \label{thm:barycntr:UrekHoppe}
The sequence $\bigl(\tfrac{U_n'}{n^2}\bigr)_n$ convergerges in distribution to a random variable $U'$, which is a solution of the fixpoint equation
\equn{ \label{eq:barycntr:Uprrek}
  U'\distreqq(1-V)^2U'+V^2U+V^2,
}
where $U'$, $V$ and $U$ are independent, $V\sim\beta(1,\theta)$ and $U$ obeys the in Corollary~\ref{thm:barycntr:Urek} characterized distribution. First and second moment of $U'$ are 
\equ{
  \Ex U'=\frac{2}{1+\theta},\qquad \Ex U'^2=\frac{12\theta+76}{3(1+\theta)(2+\theta)(3+\theta)}.
}
Again, if $\kappa$ is the normal kernel for $\sigma^2>0$, the barycentre has a mixed centered normal distribution, whose variance has expectation $\tfrac{2}{1+\theta}$ and variance $\frac{28\theta+4}{3(1+\theta)^2(2+\theta)(3+\theta)}\sigma^2$.
}

\bem{
\begin{enumerate}
  \item Basically equation~\eqref{eq:barycntr:Urek} and \eqref{eq:barycntr:Uprrek} seem to differ just by two primes, but in the second equation, $U$ is a random variable with a fixed distribution, although itself given by another fixpoint equation. Both results follow directly from limit theorems for the joint distribution of the length and the Wiener index $(T_n,W_n)$ of the trees, which is the following section dedicated to.
  \item As $\theta$ grows, expectation and variance of $U'$ tend to 0. Considering the recursive set, as $\theta$ grows, $X_0$ will be chosen as a parent more and more likely and one reaches somehow a law of large numbers regime. On the other hand, if $\theta$ vanishes, $X_0$ looses this role in favour of $X_1$ and the recursive set behaves like the one for $\theta=1$ which is randomly shifted by $X_1$. Letting formally $\theta=0$, the expectation of $U'$ is increased by by 1 compared to $U$, but both variances agree.
\end{enumerate}
}

\section{Limit theorems for length and Wiener index of Hoppe trees\label{sect:limits}}

Corollaries~\ref{thm:barycntr:Urek}~and~\ref{thm:barycntr:UrekHoppe} rely on a limit theorem for the joint distribution of $T_n$ and $W_n$ given in~\cite{rN02} for the random recursive tree, which we briefly recall: Let
\equ{
  \tilde{W}_n=\frac{W_n-\alpha_n}{n^2}, \qquad \tilde{T}_n=\frac{T_n-\gamma_n}{n},
}
be the standardized Wiener index and length of the tree. The precise values and the asymptotics of $\alpha_n$ and $\gamma_n$ are collected in Proposition~\ref{thm:barycntr:expect} with $\theta=1$. Furthermore define $\Esig(v)=v\log v+(1-v)\log(1-v)$. 

\satzn{Neininger~\cite{rN02}}{ \label{thm:gauss:gws:Neininger}
In the random recursive tree, the normalized vector $(\tilde{W}_n,\tilde{T}_n)$ converges in law together with the second moments
\equ{
  (\tilde{W}_n,\tilde{T}_n)\to (\tilde{W},\tilde{T}),
}
where the distribution of $(\tilde{W},\tilde{T})$ is the unique fixpoint of the mapping $T:\M_2\to\M_2$,
\equn{ \label{eq:gauss:gws:neininger}
  \begin{pmatrix} \tilde{W}\\ \tilde{T}\end{pmatrix} 
	\distreqq \begin{pmatrix} (1-V)^2 & V(1-V)\\ 0 & 1-V \end{pmatrix} \begin{pmatrix} \tilde{W}^\ast\\ \tilde{T}^\ast\end{pmatrix} 
	+ \begin{pmatrix} V^2 & V(1-V)\\ 0 & V \end{pmatrix} \begin{pmatrix} \tilde{W}\\ \tilde{T}\end{pmatrix}
	+ b^\ast
}
with 
\equ{
  b^\ast = \begin{pmatrix} 3V(1-V)+\Esig(V) \\ V+\Esig(V) \end{pmatrix}.
}
$(\tilde{W}, \tilde{T})$, $(\tilde{W}^\ast, \tilde{T}^\ast)$ and $V$ are independent, the two vectors equal in distribution and $V\sim \mathcal{U}[0,1]$.
}

\bem{
The structure captures the recursive structure of the tree: If we denote by $A^\ast$ the first and by $B^\ast$ the second of the two in~\eqref{eq:gauss:gws:neininger} defined matrices, then $(\tilde{W}^\ast, \tilde{T}^\ast)$ are the terms of the cutted tree containing the root and $(\tilde{W}, \tilde{T})$ belong to the tree with the root cut off.
}

Analogously one obtains a limit theorem for the Hoppe tree. Since the cut-off produces a random recursive tree and a Hoppe tree; what happens is that the Hoppe tree roughly is a pertubation of the random recursive tree. The terms of the random recursive tree take the role of a inhomogenity. Again, by $\tilde{W}_n'$ and $\tilde{T}_n'$ we denote the standardized random variables for length and Wiener index of the Hoppe tree.

\satzn{Limits in Hoppe tree}{ \label{thm:gauss:gws:Hoppe}
Let $\theta>0$. The vector $(\tilde{W}_n',\tilde{T}_n')^t$ converges in distribution and with second moments to a vector $(\tilde{W}',\tilde{T}')^t$, where its distribution is given by the fixpoint equation
\equn{ \label{eq:gauss:gws:chardistr}
  \begin{pmatrix} \tilde{W}'\\ \tilde{T}'\end{pmatrix} 
	\distreqq \begin{pmatrix} (1-V)^2 & V(1-V)\\ 0 & 1-V \end{pmatrix} \begin{pmatrix} \tilde{W}'\\ \tilde{T}'\end{pmatrix}
	+ c^\ast,
}
where
\eqa{
  c^\ast &= \begin{pmatrix} V^2 & V(1-V)\\ 0 & V \end{pmatrix} \begin{pmatrix} \tilde{W}\\ \tilde{T}\end{pmatrix}
	+ b^\ast\\
  b^\ast &= \begin{pmatrix}\left(\tfrac{\theta+5}{\theta+1}-\tfrac{2\theta+4}{\theta+1}V\right)V+\Esig(V) +\bigl(\Psi(\theta+1)- \Psi(2)\bigr) V\\
		V+\Esig(V) +\bigl(\Psi(\theta+1) -\Psi(2)\bigr) V \end{pmatrix}.
}
Here $V\sim\beta(1,\theta)^t$ and $(\tilde{W},\tilde{T})^t$ is independent of $V$ and $(\tilde{W}',\tilde{T}')^t$ and has the distribution characterized in Theorem~\ref{thm:gauss:gws:Neininger}.
}

Compared to Theorem~\ref{thm:gauss:gws:Neininger}, the inhomogenity consists of a vector $b^\ast$ which is slightly modificated, plus the limiting distribution of the vector belonging to the random recursive tree. Connecting these two results with Lemma~\ref{thm:rwrep:rec:Un}, we prove the two corollaries about the asymptotic behaviour of $U_n$ in the previous section.

\begin{proof}[Proof of Corollary~\ref{thm:barycntr:Urek}]
Denoting $Q=(-1,1)$, we have
\equn{ \label{eq:gauss:gws:Un_norm}
  \frac{U_n}{n^2}=Q{\tilde{W}_n\choose\tilde{T}_n}+\frac{n\Ex T_n-\Ex W_n}{n^2},
}
where $\tilde{W}_n=\tfrac{W_n-\Ex W_n}{n^2}$ and $\tilde{T}_n=\tfrac{T_n-\Ex T_n}{n}$. Since the second term converges as well as the vector $(\tilde{W}_n,\tilde{T}_n)^t$ converges in distribution and with the first two moments, also $\tfrac{U_n}{n^2}$ converges to some $U=Q (\tilde{W},\tilde{T})^t+1$ in distribution with
\eqa{
  U &\distreqq Q\bigl[ B^\ast{W\choose T}+A^\ast{W^\ast\choose T^\ast}+b^\ast\bigr]+1\\
	&= V^2 U^\ast+(1-V)^2U-V^2-(1-V)^2+1+V-3V(1-V)\\
	&= V^2 U^\ast+(1-V)^2U+V^2
}
because of
\eqa{
  QB^\ast &=(-V^2,-V(1-V)+V)=V^2Q\\
  QA^\ast &=(-(1-V)^2,-V(1-V)+(1-V))=(1-V)^2Q\\
  Qb^\ast &= -3V(1-V)-\Esig(V)+V+\Esig(V)=V-3V(1-V).
}

By using the distributional equality and independence, we get $\Ex U=1$ and $\Ex U^2=\tfrac{11}{9}$.
\end{proof}

An analogue reasoning yields the convergence for the Hoppe tree.

\begin{proof}[Proof of Corollary~\ref{thm:barycntr:UrekHoppe}]
The ansatz~\eqref{eq:gauss:gws:Un_norm} still holds with
\equ{
  U'=T'-W'+\frac{2}{1+\theta}.
}
Therefore,
\eqa{
  U'&\distreqq QA^\ast{\tilde{W}'\choose \tilde{T}'}+QB^\ast{\tilde{W}\choose \tilde{T}}+Qb^\ast+\frac{2}{1+\theta}\\
	&=(1-V)^2Q{\tilde{W}'\choose \tilde{T}'}+V^2Q{\tilde{W}\choose \tilde{T}}-\left[\frac{4}{\theta+1}-\frac{2\theta+4}{\theta+1}V\right]V+\frac{2}{1+\theta}\\
	&=(1-V)^2U'+V^2U-\left[\frac{4}{\theta+1}-\frac{2\theta+4}{\theta+1}V\right]V+\frac{2}{1+\theta}-\frac{2}{1+\theta}(1-V)^2-V^2\\
	&=(1-V)^2U'+V^2U+V^2
}
To obtain the moments, observe that
\equ{
  \Ex U'=\frac{\theta}{2+\theta}\Ex U'+\frac{4}{(\theta+1)(\theta+2)},
}
i.e. $\Ex U'=\tfrac{2}{1+\theta}$ and moreover, 
\equ{
  \Ex U'^2=\Ex U'^2\frac{\theta\Gamma(\theta+4)}{\Gamma(\theta+5)}+\frac{16\theta\Gamma(\theta+2)}{(\theta+1)\Gamma(\theta+5)}+\frac{4!\Gamma(\theta+1)}{\Gamma(\theta+5)}\frac{38}{11}
}
i.e. $\Ex U'^2=\tfrac{4\theta+\frac{76}{3}}{(1+\theta)(2+\theta)(3+\theta)}$.
\end{proof}

Before we turn to the proof of Theorem~\ref{thm:gauss:gws:Hoppe}, a few remaks on Proposition~\ref{thm:barycntr:expect}.

\begin{proof}[Proof of Proposition~\ref{thm:barycntr:expect}]
Firstly remark that by conditioning,
\equ{
  \Ex[D_{in-1}|\mathcal{J}_{n-1}] =\frac{\theta}{\theta+n-2}D_{i0} 
  		+\frac{1}{\theta+n-2}\sum_{j=1}^{n-2}D_{ij}
  		=\frac{2}{\theta+n-2}R_{n-1}+\frac{1}{\theta+n-2}T_{n-1}.
}
$D_{i0}=0$ since the last common ancestor is $\l_{i0}=0$. If $n-1$  grows at $j$, then given the tree the last common ancestor of $i$ and $n-1$ is exactly the one of $i$ and $j$. Immedialty,
\eqa{
  \Ex[R_n|\mathcal{J}_{n-1}] &= R_{n-1}+\Ex[\sum_{i=0}^{n-1}D_{in-1}|\mathcal{J}_{n-1}]\\
		&=\frac{\theta+n}{\theta+n-2}R_{n-1}+\frac{1}{\theta+n-2}T_{n-1}.
}
Combining~\eqref{en:rek:Tn},~\eqref{en:rek:Rn} and Lemma~\ref{thm:rwrep:rec:Un}, we obtain~\eqref{en:rek:Wn}.


Since the expectation of $T_n$ is known, we get $\Ex[U_n]=\tfrac{\theta+n}{\theta+n-2}\Ex[U_{n-1}]+h_{n-1}^\theta+1$. Divide by $(\theta+n)(\theta+n-1)$, then because of $\Ex[U_1]=0$,
\eqa{
  \Ex \frac{U_n}{(\theta+n)(\theta+n-1)} &= \sum_{j=1}^{n-1}\frac{1}{(\theta+j)(\theta+j+1)}(1+h_{j}^\theta)\\
  	&= \sum_{j=1}^{n-1}\frac{1}{(\theta+j)(\theta+j+1)}
  			+\sum_{j=1}^{n-1}\sum_{k=1}^{j-1}\frac{1}{(\theta+j)(\theta+j+1)}\frac{1}{\theta+k}\\
  	&= \frac{1}{\theta+1}-\frac{1}{\theta+n}
  			+\sum_{k=1}^{n-2}\frac{1}{\theta+k}\sum_{j=k+1}^{n-1}\frac{1}{(\theta+j)(\theta+j+1)}\\
  	&= \frac{1}{\theta+1}-\frac{1}{\theta+n}
  			+\sum_{k=1}^{n-2}\frac{1}{\theta+k}\left(\frac{1}{\theta+k+1}-\frac{1}{\theta+n}\right).
}
Collecting the terms we get the expression for $\Ex U_n$, the asymptotic behaviour follows from $\log n-h_n^\theta\to\Psi(\theta+1)$.

Analogously we obtain the expression for $\Ex W_n$ by some calculus from
\equ{
  \Ex \frac{W_n}{(\theta+n)(\theta+n-1)} = \sum_{j=1}^{n-1}\frac{1}{(\theta+j)(\theta+j+1)}\bigl[(\theta-1)h_{j}^\theta+j\bigr].
}
\end{proof}

The proof of Theorem~\ref{thm:gauss:gws:Hoppe} basically goes along the same lines as the proof of~\ref{thm:gauss:gws:Neininger}, it is an application of the contraction method: Find a recursion for the finite trees via cutting off a branch and then show that the parameters obtained behave asymptotically in a sufficiently nice way. We need to identify the distribution of the size of the subtrees.

\lemma{
Let $R_n$ and $U_n$ be given as in Lemma~\ref{thm:rwrep:rec:Un}, then
\eqan{
  2R_n &\distreqq 2R_K+2R_{n-K}'+K(2n-K-1)\\
  U_n &\distreqq U_K+U_{n-K}'+K^2+2K(n-K),\label{eq:martingal:Urec}
}
where $R_j$ and $R_j'$, $U_j$ and $U_j'$, $j=1,\ldots,n$, are independent and independent of $K$ with distribution
\equn{ \label{eq:limits:distrK}
  \Pr(K=m)={n-2\choose m-1}\frac{\theta^{(n-m-1)}(m-1)!}{(\theta+1)^{(n-2)}},\qquad m\in\{1,\ldots,n-1\}
}
and for $j=1,\ldots,n$, the underlying distribution for $U_j$ is the random recursive tree, and for $U_j'$ the Hoppe tree.
}

\bem{
\begin{enumerate}
  \item For $\theta=1$, $K$ is uniformly distributed on $\{1,\ldots,n-1\}$, and $U_j\distreqq U_j'$. This is the result in~\cite{DF99}, Lemma 2.1.
  \item More generally, if one initializes P\'olya's urn with a white ball with weight $\theta$ and a red ball with weight 1 and adds in each step one ball with weight 1, then $K$ is the number of red balls after $n-2$ draws. 
\end{enumerate}
}
\begin{proof}
Recall that $2R_n$ is the sum of the distance of the common ancestor to the root over all ordered pairs of distinct vertices. Let $K$ be the size of the subtree originated at 1 und decompose $R_n$ into the sum over pairs in this subtree, the sum over the pairs of the subtree rooted at 0 and not containing 1, and finally the sum over pairs with one partner of each subtree. Conditional $K$, they yield $2R_K+K(K-1)$, since every pair yields an additional 1 from 1 to the root, $2R_{n-K}$ for the other branch and finally $2K(n-K)$ for all mixtures. \eqref{eq:martingal:Urec} follows analogously. 

The distribution of $K$ given by~\eqref{eq:limits:distrK} follows e.g. inductively from the binomial theorem for rising factorials.
\end{proof}

Now we collected the tools to prove the limit theorem for length and Wiener index for the Hoppe tree.

\begin{proof}[Proof of theorem~\ref{thm:gauss:gws:Hoppe}]
For $W_n'$ the fundamental recursion is
\equ{
  W_n'\distreqq W_{n-I_n}'+W_{I_n}+I_n T_{n-I_n}'+(n-I_n)T_{I_n}+I_n(n-I_n).
}
Hence
\equ{
  \begin{pmatrix} W_n'\\ T_n'\end{pmatrix} 
	\distreqq \begin{pmatrix} 1 &I_n\\ 0 & 1\end{pmatrix} \begin{pmatrix} W_{n-I_n}'\\ T_{n-I_n}'\end{pmatrix}
	+ \begin{pmatrix} 1 &n-I_n\\ 0 & 1\end{pmatrix} \begin{pmatrix} W_{I_n}\\ T_{I_n}\end{pmatrix}
	+ \begin{pmatrix} I_n(n-I_n)\\ I_n\end{pmatrix}.
}
Moreover, denote the expectations of Wiener index and total length by $\alpha_n'=\Ex[W_n']$ and $\gamma_n'=\Ex[T_n']$. Then normalization yields
\equ{
  \begin{pmatrix} \tilde{W}_n'\\ \tilde{T}_n'\end{pmatrix} 
	\distreqq A_n^\ast \begin{pmatrix} \tilde{W}_{n-I_n}'\\ \tilde{T}_{n-I_n}'\end{pmatrix}
	+ B_n^\ast \begin{pmatrix} \tilde{W}_{I_n}\\ \tilde{T}_{I_n}\end{pmatrix}
	+ b_n^\ast,
}
where
\eqa{
  A_n^\ast &= \begin{pmatrix} \tfrac{1}{n^2} & 0 \\ 0 & \tfrac{1}{n}\end{pmatrix} \begin{pmatrix} 1 &I_n\\ 0 & 1\end{pmatrix}\begin{pmatrix} (n-I_n)^2 &I_n(n-I_n)\\ 0 & n-I_n\end{pmatrix} =\begin{pmatrix} \Bigl[1-\tfrac{I_n}{n}\Bigr]^2 & \tfrac{I_n}{n}\Bigl[1-\tfrac{I_n}{n}\Bigr]\\ 0 & 1- \tfrac{I_n}{n}\end{pmatrix}\\
  B_n^\ast &= \begin{pmatrix} \tfrac{1}{n^2} &0\\ 0 & \tfrac{1}{n}\end{pmatrix} \begin{pmatrix} 1 & n-I_n\\ 0 & 1\end{pmatrix}\begin{pmatrix} I_n^2 & 0 \\ 0 & I_n \end{pmatrix} =\begin{pmatrix} \Bigl[\tfrac{I_n}{n}\Bigr]^2 & \tfrac{I_n}{n}\Bigl[1-\tfrac{I_n}{n}\Bigr]\\ 0 & \tfrac{I_n}{n}\end{pmatrix}\\
  b_n^\ast &= \begin{pmatrix} I_n(n-I_n)\\ I_n\end{pmatrix}-\begin{pmatrix} \tfrac{\alpha_n'}{n^2} \\ \tfrac{\gamma_n'}{n} \end{pmatrix}+ A_n^\ast \begin{pmatrix} \tfrac{\alpha_{n-I_n}'}{(n-I_n)^2} \\ \tfrac{\gamma_{n-I_n}'}{n-I_n} \end{pmatrix} + B_n^\ast \begin{pmatrix} \tfrac{\alpha_{I_n}}{I_n^2} \\ \tfrac{\gamma_{I_n}}{I_n} \end{pmatrix}
}
Thus we identified the recursion and now have to show that $(A_n^\ast)_{n}$ as well as the inhomogenity converge suitably. First of all,
\eqa{
  b_{n1}^\ast &= \begin{multlined}[t] 
		\tfrac{I_n}{n}\Bigl[1-\tfrac{I_n}{n}\Bigr] +\Bigl[1-\tfrac{I_n}{n}\Bigr]^2\log\Bigl[1-\tfrac{I_n}{n}\Bigr] +\tfrac{I_n}{n}\Bigl[1-\tfrac{I_n}{n}\Bigr]\log\Bigl[1-\tfrac{I_n}{n}\Bigr] +\Bigl[\tfrac{I_n}{n}\Bigr]^2\log\Bigl[\tfrac{I_n}{n}\Bigr]\\
		+\tfrac{I_n}{n}\Bigl[1-\tfrac{I_n}{n}\Bigr]\log\Bigl[\tfrac{I_n}{n}\Bigr] + \tfrac{2}{\theta+1}\biggl(1-\Bigl[1-\tfrac{I_n}{n}\Bigr]^2\biggr)\\
		+\Psi(\theta+1) \biggl(1-\Bigl[1-\tfrac{I_n}{n}\Bigr]^2-\tfrac{I_n}{n}\Bigl[1-\tfrac{I_n}{n}\Bigr]\biggr) -\psi(2)\biggl(\tfrac{I_n}{n}\Bigl[1-\tfrac{I_n}{n}\Bigr]+\Bigl[\tfrac{I_n}{n}\Bigr]^2\biggr)-\Bigl[\tfrac{I_n}{n}\Bigr]^2 +o(1)
	\end{multlined}	\\
 &\to \begin{multlined}[t]
		V(1-V)+(1-V)^2\log(1-V) +V(1-V)\log(1-V)+V^2\log V+ V(1-V)\log V\\
		+ \tfrac{2}{\theta+1}V(2-V) +\Psi(\theta+1)V- \Psi(2) V -V^2
	\end{multlined} \\
  b_{n2}^\ast &= \begin{multlined}[t]
		\tfrac{I_n}{n}+\Bigl[1-\tfrac{I_n}{n}\Bigr]\log\Bigl[1-\tfrac{I_n}{n}\Bigr]+\tfrac{I_n}{n}\log\tfrac{I_n}{n} +\Psi(\theta+1) \biggl(1-\Bigl[1-\tfrac{I_n}{n}\Bigr]\biggr)-\Psi(2) \tfrac{I_n}{n}+o(1)
	\end{multlined}\\
 &\to V+(1-V)\log(1-V)+V\log V +\Psi(\theta+1)V -\Psi(2) V
}
The $o(1)$ terms are random but vanish uniformly. Since we can choose a probability space such that $\tfrac{I_n}{n}$ converges almost surely and thus in $L^2$ to a $\beta(1,\theta)$ distributed random variable, we obtain by dominated convergence the convergence of $b_n^\ast$ in $L^2$ to $b^\ast$ as well as of $A_n^\ast$ and $B_n^\ast$ to $A^\ast$ and $B^\ast$. Remark that by Theorem~\ref{thm:gauss:gws:Neininger}, the joint probability space can be adjusted such that $(\tilde{W}_{n},\tilde{T}_{n})$ also converges $L^2$.

Secondly, denoting by $\lVert\,\cdot\,\rVert$ the spectral norm, $\Ex[\lVert(A^\ast)^{t}A^\ast\rVert]<1$ holds for all $\theta>0$ since the eigenvalues of $(A^\ast)^{t}A^\ast$ are
\equ{
  \lambda(V)=(1-V)^2\left\{1+V^2-V\left(1\pm\sqrt{(1-V)^2+1}\right)\right\},
}
the larger one in absolute value being the one with $-$. But then, by estimating the root by $2$,
\equ{
  \Ex[\lVert(A^\ast)^{t}A^\ast\rVert]=\Ex[\lambda(V)]\leq \frac{\theta}{2+\theta}\left[1+\frac{1}{3+\theta} + \frac{2}{(4+\theta)(3+\theta)}\right],
}
which is less than 1 for all $\theta>0$ as a short comparison shows.

Finally $\Ex[1_{\{n-I_n\leq l\}\cup\{n-I_n=n\}}\lVert(A_n^\ast)^{t}A_n^\ast\rVert]$ vanishes for all fixed $l\in\N$ since the eigenvalues of $(A_n^\ast)^{t}A_n^\ast$ are bounded by $3$ and for fixed $l$, 
\equ{
  \Pr(n-I_n\leq l,n-I_n=n)\to 0.
}

Thus~\cite[Theorem 4.1]{rN01} applies und $(\tilde{T_n},\tilde{W_n})^t$ converges in the given sense to some $(\tilde{T},\tilde{W})^t$, whose distribution is characterized by equation~\eqref{eq:gauss:gws:chardistr}.


\end{proof}




\bibliographystyle{alpha-abbrv}


\begin{thebibliography}{10}

\bibitem{BZ11}
Bach, A. and Zessin, H. (2011).
\newblock The particle structure of the quantum mechanical Bose and Fermi gas.
\newblock {\em preprint}.

\bibitem{DF99}
Dobrow, R.~P. and Fill, J.~A. (1999).
\newblock Total Path Length for Random Recursive Trees.
\newblock {\em Combin. Probab. Comput.}, 8:317--33.

\bibitem{fH87}
Hoppe, F.~M. (1987).
\newblock The sampling theory of neutral alleles and an urn model in population genetics,
\newblock {\em J. Math. Biol.}, 25(2):123--59.

\bibitem{fH84}
Hoppe, F.~M. (1984).
\newblock Pólya-like urns and the Ewens' sampling formula.
\newblock {\em J. Math. Biol.}, 20(1):91--4.

\bibitem{oK78}
Kallenberg, O. (1978).
\newblock On conditional intensities of point processes.
\newblock {\em Z. Wahrscheinlichkeitstheorie verw. Geb.}, 41:205--20.

\bibitem{MKM78}
Kerstan, J., Matthes, K. and Mecke, J. (1978).
\newblock {Infinitely Divisible Point Processes}.
\newblock John Wiley \& Sons.

\bibitem{LN12}
Leckey, K. and Neininger, R. (2012).
\newblock Asymptotic analysis of Hoppe trees.
\newblock {\em Preprint}.

\bibitem{hM91}
Mahmoud, H.~M. (1991).
\newblock Limiting distributions for path lengths in recursive trees.
\newblock {\em Probab. Engrg. Inform. Sci.}, 5(1):53--9.

\bibitem{RN12}
Nehring, B. and Rafler, M. (2012).
\newblock The branching Pólya sum process.
\newblock {\em Preprint}.

\bibitem{rN01}
Neininger, R. (2001).
\newblock On a multivariate contraction method for random recursive structures with applications to quicksort.
\newblock {\em Random Structures and Algorithms}, 19(3--4):498--524.

\bibitem{rN02}
Neininger, R. (2002).
\newblock The Wiener index of random trees.
\newblock {\em Combin. Probab. Comput.}, 11(6):587--97.

\bibitem{hZ09}
Zessin, H. (2009).
\newblock Der Papangelou Prozess.
\newblock {\em Journal of Contemporary Mathematical Analysis}, 44(1):36--44.

\end{thebibliography}

\end{document}